\documentclass [11pt]{amsart}
\usepackage{array}
\usepackage{amsmath,amssymb,graphics,graphicx}

\date{} 

\newtheorem{theorem}{Theorem}[section]

\newtheorem{lemma}[theorem]{Lemma}

\newtheorem{corollary}[theorem]{Corollary}

\newcommand{\ignore}[1]{}
\newcommand{\tr}{{\operatorname{tr}}}
\newcommand{\monic}{{\operatorname{monic}}}
\newcommand{\pr}{{\operatorname{prime}}}
\newcommand{\dis}{{\operatorname{distinct}}}
\newcommand{\lalin}{{\operatorname{i}}}
\newcommand{\mob}{{\operatorname{o}}}

\begin{document}
\title[ Products of traces of high powers of the Frobenius class ]
{Statistics for products of traces of high powers of the Frobenius class of hyperelliptic
curves}
\author{Edva Roditty-Gershon}
\address {Raymond and Beverly Sackler School of Mathematical Sciences,
Tel Aviv University, Tel Aviv 69978, Israel}
\email{rodittye@post.tau.ac.il}
\date{\today} 
\thanks{Supported in part by the Israel Science Foundation (grant No. 1083/10).}
\maketitle

\maketitle
\begin{abstract}
We study the averages of products of traces of high powers of the
Frobenius class of hyperelliptic curves of genus g over a fixed
finite field. We show that for increasing genus g, the limiting
expectation of these products equals to the expectation when the
curve varies over the unitary symplectic group USp(2g). We also
consider the scaling limit of linear statistics for eigenphases of
the Frobenius class of hyperelliptic curves, and show that their
first few moments are Gaussian.
\end{abstract}
\tableofcontents
\clearpage
\setcounter{page}{1}
\section{Introduction}
Let $\emph{C}$ be a nonsingular projective curve of genus $g$, defined over a finite field $\mathbb{F}_{q}$ of odd cardinality $q$. The zeta function of $\emph{C}$ is defined as
$$Z_{\emph{C}}(u)=\exp\sum_{n=1}^{\infty}N_{n}(\emph{C})\frac{u^{n}}{n},~~~~|u|<\frac{1}{q}\eqno(1.1)$$
where $N_{n}(\emph{C})$ is the number of points of $\emph{C}$ with coefficients in an extension $\mathbb{F}_{q^{n}}$ of $\mathbb{F}_{q}$ of degree n.
The zeta function is known to be a rational function of the form
$$Z_{\emph{C}}(u)=\frac{P_{\emph{C}}(u)}{(1-u)(1-qu)}\eqno(1.2)$$
where $P_{\emph{C}}\in \mathbb{Z}[u]$ is a polynomial of degree $2g$, with $P_{\emph{C}}(0)=1$ , satisfying a functional equation
$$P_{\emph{C}}(u)=(qu^{2})^{g}P_{\emph{C}}(\frac{1}{qu}).$$
By the Riemann hypothesis (proved by Weil \cite{weil}) ,We may interpret $P_{\emph{C}}(u)$ as the characteristic polynomial of a $2g\times 2g$ unitary matrix $\Theta_{\emph{C}}$, Where the eigenvalues $e^{i\theta_{j}}$ of $\Theta_{\emph{C}}$ correspond to zeros $q^{-1/2}e^{-i\theta_{j}}$ of $Z_{\emph{C}}(u)$:
$$P_{\emph{C}}(u)=\det(I-u\sqrt{q}\Theta_{\emph{C}})\eqno(1.3)$$
The conjugacy class of $\Theta_{\emph{C}}$ is called the $\emph{unitarized Frobenius class}$ of $\emph{C}.$
\\We consider the family $\mathcal{H}_{2g+1}$ of hyperelliptic curves of genus $g$ given in affine
form by an equation
$$\emph{C}_{Q}: y^{2}=Q(x)$$
where $Q(x)\in\mathbb{F}_{q}[x]$ is a square-free, monic polynomial
of degree $2g+1$. We will study the expected value of products of
traces of high powers of the Frobenius class of $\emph{C}$ as we
vary the curve $\emph{C}$ over $\mathcal{H}_{2g+1},$ and show that
these statistics determine the n-level density of the eigenvalues.
Our work is in the limit of large genus and fixed constant field.
\\Consider $\mathcal{H}_{2g+1}$ as a probability space with the uniform probability measure, so that the expected value of any function $F$ on $\mathcal{H}_{2g+1}$ is defined as
$$\langle F \rangle:=\frac{1}{\# \mathcal{H}_{2g+1}}\sum_{Q\in \mathcal{H}_{2g+1}}F(Q).$$

Katz and Sarnak \cite{k-s} showed that for a fixed genus, the Frobenius classes $\Theta_{Q}$ become uniformly distributed in $USp(2g)$ in the limit $q\rightarrow\infty$ of large field size. That is, for any continuous function on the space of conjugacy classes of $USp(2g)$,
$$\lim_{q\rightarrow\infty}\langle F(\Theta_{Q})\rangle=\int_{USp(2g)}F(U)dU.$$

If we take the opposite limit, that of fixed constant field and large genus
$g\rightarrow\infty$ (that is without first taking $q\rightarrow\infty$, which was crucial to the approach of Katz
and Sarnak), since the matrices $\Theta_{Q}$ now inhabit different spaces as $g$ grows, it is not clear how to formulate an equidistribution problem. However, we can discuss the statistics of products of traces of powers of $\Theta_{Q},$ that is, $\langle \prod_{j=1}^{n} (\tr U^{k_{j}})^{a_{j}}\rangle.$
Rudnick \cite{rud} showed that for a fixed constant field and large genus $g\rightarrow\infty$, if $3\log_{q}g<n<4g-5\log_{q}g$ but $n\neq 2g$ then
$$\langle \tr U^{n}\rangle=\int_{USp(2g)}\tr U^{n}dU+o(\frac{1}{g}).$$
In the case of fixed $k_{1},\ldots,k_{n},a_{1},\ldots,a_{n}$ Bucur-David-Feigon-Lal$\acute{\lalin}$n \cite{Bu-Da-Feigon-Lal} studied the variation of the trace of the Frobenius endomorphism in the cyclic trigonal
ensemble. They showed that for $q$ fixed and $g$ increasing, the limiting distribution of the trace of Frobenius equals the sum of
q +1 independent random variables taking the value $0$ with probability $2/(q +2)$ and $1, e^{2\pi i/3}, e^{4\pi i/3}$
each with probability $q/(3(q + 2))$. This extends the work of Kurlberg and Rudnick \cite{k-r} who considered
the same limit for hyperelliptic curves.

In this paper (in continuation of Rudnick's work \cite{rud}), we study the general case of average of product of traces $\langle \prod_{j=1}^{n} (\tr U^{k_{j}})^{a_{j}}\rangle,$ where $k_{1},\ldots,k_{n}$ are of order of the genus $g$, and $a_{1},\ldots,a_{n}$ are fixed .

\subsection{Result}
First we state a result \cite{d-s},\cite{d-e},\cite{rud-hug} which
expresses the mean value of products of traces of high powers when
averaged over the unitary symplectic group $USp(2g)$ in terms of
independent standard normal random variables.
\\
\\Let $Z_{j}$ be independent standard normal random variables, and let
$$\eta_{k_{j}}=\left\{ \begin{array}{ll}
1 &\mbox{if} ~k_{j}~ \mbox{is even}\\
0 &\mbox{if} ~k_{j}~ \mbox{is odd}
\end{array}
\right.$$
If $k_{j},a_{j}\in \{1,2,\ldots\}$ for $1\leq j \leq n$ are such that $\sum_{j=1}^{n}a_{j}k_{j}\leq 2g+1$, $k_{j}$ distinct, then
$$\int_{{USp(2g)}} \prod_{j=1}^{n}(\tr U^{k_{j}})^{a_{j}}dU=\mathbb{E} (\prod_{j=1}^{n}(\sqrt{k_{j}}Z_{j}-\eta_{k_{j}})^{a_{j}})$$
Where $\mathbb{E}$ denotes the expectation.
For the proof see \cite{d-s},\cite{d-e},\cite{rud-hug}.
\\Since $Z_{j}$ are independent standard normal random variables, we have
$$\mathbb{E}(\prod_{j=1}^{n}(\sqrt{k_{j}}Z_{j}-\eta_{k_{j}})^{a_{j}})=\prod_{j=1}^{n}\mathbb{E}( (\sqrt{k_{j}}Z_{j}-\eta_{k_{j}})^{a_{j}})=\prod_{j=1}^{n}\mathbb{E}( \sum_{i=0}^{a_{j}}{a_{j}\choose i}(\sqrt{k_{j}}Z_{j})^{i}(-\eta_{k_{j}})^{a_{j}-i})=$$
$$\prod_{j=1}^{n} \sum_{i=0}^{a_{j}}{a_{j}\choose i}(\sqrt{k_{j}})^{i}\mathbb{E}((Z_{j})^{i})(-\eta_{k_{j}})^{a_{j}-i}=
\prod_{j=1}^{n} \sum_{i=0}^{\lfloor\frac{a_{j}}{2}\rfloor}{a_{j}\choose 2i}(k_{j})^{i}(\frac{(2i)!}{2^{i}(i)!})(-\eta_{k_{j}})^{a_{j}-2i}\eqno(1.4)$$
\\We will show
\begin{theorem}
Assume $k_{j}\in \{1,2,\ldots\}$ for $1\leq j \leq n$ are such that $\sum_{j=1}^{n}a_{j}k_{j}\leq 2g-1$ for fixed integers $a_{j}$. Assume that $k_{j}$ are distinct and $\log_{q}g\ll\min(k_{1},\ldots,k_{n})$, then
$$\langle \prod_{j=1}^{n} (\tr U^{k_{j}})^{a_{j}}\rangle=
\prod_{j=1}^{n}\sum_{i_{j}=0}^{\lfloor\frac{a_{j}}{2}\rfloor}(k_{j})^{i_{j}}{a_{j}\choose
2i_{j}}\frac{(2i_{j})!}{2^{i_{j}}(i_{j})!}(-\eta_{k_{j}})^{a_{j}-2i_{j}}
+o(1)\eqno(1.5)$$
\end{theorem}
Comparing (1.4) and (1.5) we find
\begin{corollary}
If $\log_{q}g\ll\min(k_{1},\ldots,k_{n})$ and $\sum_{j=1}^{n}k_{j}a_{j}\leq2g-1$, then
$$\langle \prod_{j=1}^{n}(\tr U^{k_{j}})^{a_{j}}\rangle= \int_{USp(2g)}\prod_{j=1}^{n}(\tr U^{k_{j}})^{a_{j}}dU+o(1)\eqno(1.6)$$
\end{corollary}

to prove these results, we can not use the same methods that were used for the fixed genus case by Katz and Sarnak \cite{k-s}. Rather, we use a variant of the analytic methods similar to those used in \cite{rud}.
\subsection{Application: The n-level density}
Denote by $\theta_{1},\ldots,\theta_{N}$ the sequence of angles of
$U$ a unitary matrix of size $N\times N.$ The traces of powers
determine the number of sets of angles
$\theta_{i_{1}},\ldots,\theta_{i_{n}}$ lying in a subinterval of
$\mathbb{R}/2\pi \mathbb{Z}$, or the $n$-level density. For the case
of $n=1$ or the one-level density see \cite{rud}. To define the
$n$-level density, we start with an even test function $f$, in the
Schwartz space $\mathcal{S}(\mathbb{R}),$ and for any $N\geq 1$ set
$$F(\theta):=\sum_{k\in \mathbb{Z}}f(N(\frac{\theta}{2\pi}-k)),$$
which has a period of $2\pi$ and is localized in an interval of size $\approx 1/N$ in $\mathbb{R}/2\pi \mathbb{Z}.$ For a unitary $N\times N$ matrix $U$ with eigenvalues $e^{i\theta_{j}}$ , $j=1,\ldots,N,$ define
$$Z_{f}(U):=\sum_{j=1}^{N}F(\theta_{j}),$$
which counts the number of "low-lying" eigenphases $\theta_{j}$ in the smooth interval of length $\approx 1/N$ around the origin defined by $f.$
The product $Z_{f}^{n}$ counts the number of sets of angles $\theta_{i_{1}},\ldots,\theta_{i_{n}}$ in the smooth interval of length $\approx 1/N$ around the origin defined by $f.$ In order to study the $n$-level density, we need to compute the $n^{th}$ moment of $Z_{f}.$

Katz and Sarnak \cite{k-s} conjectured that for fixed $q$, the expected value of $Z_{f}$ over $\mathcal{H}_{2g+1}$ will converge to $\int_{USp(2g)}Z_{f}(U)dU$ as $g\rightarrow\infty$ for any such test function $f.$ Rudnick \cite{rud} proved this conjecture for a test function $f$ such that the Fourier transform $\widehat{f}$ supported in $(-2,2).$
Corollary 1.2 implies:
\begin{corollary}
If supp$\hat{f}\subseteq(\frac{-1}{m},\frac{1}{m})$ then the first m moments of $Z_{f}(U)$ converge to the Gaussian moments with mean
$$\hat{f}(0)-\int_{0}^{1}\hat{f}(u)du$$
and variance
$$2\int_{-1/2}^{1/2}|u|\hat{f}(u)^{2}du.$$
\end{corollary}
This is called "Mock Gaussian" behavior in \cite{rud-hug}.
\\To show Corollary 1.3, one uses a Fourier expansion to see that (for $N=2g$)
$$Z_{f}(U)=\int_{-\infty}^{\infty}f(x)dx+\frac{1}{N}\sum_{k\neq 0}\widehat{f}(\frac{k}{N})\tr U^{k}\eqno(1.7)$$
and then by Corollary 1.2 and  \cite{rud-hug}, the above follows.
\section{Background on Dirichlet characters and L-functions}
In this section we review some known background on Quadratic L-function. See \cite{ros} for details.
\subsection{The zeta function} For a nonzero polynomial $f\in \mathbb{F}_{q}[x]$, we define the norm $|f|:=q^{\deg f}.$ A prime polynomial is a monic irreducible polynomial.
For a monic polynomial $f$, the von Mangoldt function $\Lambda(f)$ is defined to be zero unless $f$ is a prime power in which case $\Lambda(P^{k})=\deg P.$

The analog of Riemann's zeta function is
$$\zeta_{q}(s):=\prod_{P~\pr}(1-|P|^{-s})^{-1},~~~~\mathcal{R}(s)<1\eqno(2.1)$$
As a result of expanding in additive form using unique factorization, we have
$$\zeta_{q}(s)= \frac{1}{1-q^{1-s}}\eqno(2.2)$$
The following identity is equivalent to (2.2):
$$\sum_{\substack{\deg f=n\\ f ~\monic}}\Lambda(f)=q^{n}\eqno(2.3)$$

Let $\pi_{q}(n)$ be the number of prime polynomials of degree $n$. The Prime Polynomial Theorem in $\mathbb{F}_{q}[x]$ asserts that
$$\pi_{q}(n)=\frac{q^{n}}{n}+O(q^{n/2})\eqno(2.4)$$
which follows from (2.3).
\subsection{Quadratic characters} Let $P\in\mathbb{F}_{q}[x]$ ($q$ odd) be a prime polynomial. The quadratic
residue symbol $(\frac{f}{P})\in \{\pm1\}$ is defined for $f$ coprime to $P$ by
$$(\frac{f}{P})\equiv f^{(\frac{|P|-1}{2})} \mod P$$
For arbitrary monic $Q\in \mathbb{F}_{q}[x]$ and for $f$ coprime to $Q$, the Jacobi symbol $(\frac{f}{Q})$ is defined by writing
$Q=\prod P_{j}$ as a product of prime polynomials and setting
$$(\frac{f}{Q})=\prod(\frac{f}{P_{j}})$$
If $f,Q$ are not coprime we set $(\frac{f}{Q})= 0$.

The law of quadratic reciprocity asserts that for $A,B\in\mathbb{F}_{q}[x]$ monic polynomials
$$(\frac{B}{A})=(-1)^{(\frac{q-1}{2})\deg A\deg B}(\frac{A}{B})\eqno(2.5)$$
For $D\in \mathbb{F}_{q}[x]$ a monic polynomial of positive degree which is not a perfect square, we define the quadratic character $\chi_{D}$ by
$$\chi_{D}=(\frac{D}{f})\eqno(2.6)$$
\subsection{L-functions} For the quadratic character $\chi_{D}$, the corresponding L-function is defined by
$$\mathcal{L}(u,\chi_{D}):=\prod_{P~\pr}(1-\chi_{D}(P)u^{\deg P})^{-1},~~~~|u|<\frac{1}{q}.$$
Expanding in additive form using unique factorization, we write
$$\mathcal{L}(u,\chi_{D})=\sum_{\beta\geq 0}A_{D}(\beta)u^{\beta}$$
with
$$A_{D}(\beta):=\sum_{\substack{\deg B=\beta\\B~\monic}}\chi_{D}(B).$$
If $D$ is nonsquare of positive degree, then $A_{D}(\beta)=0$ for $\beta\geq\deg D$ and hence the L-function is in fact a polynomial of degree at most $\deg D-1.$

Now, assume that $D$ is also square-free. Then $\mathcal{L}(u,\chi_{D})$ has a trivial zero at $u=1$ if and only if $\deg D$ is even. Thus
$$\mathcal{L}(u,\chi_{D})=(1-u)^{\lambda}\mathcal{L}^{\ast}(u,\chi_{D})=,~~~~\lambda=\left\{ \begin{array}{ll}
1 &\mbox{$\deg D ~even$}\\
0 &\mbox{$\deg D~ odd$}
\end{array}
\right.$$
where $\mathcal{L}^{\ast}(u,\chi_{D})$ is a polynomial of even degree
$$2\delta=\deg D-1-\lambda$$
satisfying the functional equation
$$\mathcal{L}^{\ast}(u,\chi_{D})=(qu^{2})^{\delta}\mathcal{L}^{\ast}(\frac{1}{qu},\chi_{D}).$$
We write
$$\mathcal{L}^{\ast}(u,\chi_{D})=\sum_{\beta=0}^{2\delta}A^{\ast}_{D}(\beta)u^{\beta},$$
where $A^{\ast}_{D}(0)=1,$ and the coefficients $A^{\ast}_{D}(\beta)$ satisfy
$$A^{\ast}_{D}(\beta)=q^{\beta-\delta}A^{\ast}_{D}(2\delta-\beta)\eqno(2.7)$$
In particular, the leading coefficient is $A^{\ast}_{D}(2\delta)=q^{\delta}.$
\subsection{The explicit formula}
For $D$  monic, square-free, and of positive degree, the zeta function (1.2) of the hyperelliptic curve $y^{2}=D(x)$ is
$$Z_{D}(u)=\frac{\mathcal{L}^{\ast}(u,\chi_{D})}{(1-u)(1-qu)}.$$
By the Riemann Hypothesis (proved by Weil \cite{weil}) all the zeros of $Z_{D}(u)$, hence of $\mathcal{L}^{\ast}(u,\chi_{D})$, lie on the circle $|u|=1/q.$ Thus we may write
$$\mathcal{L}^{\ast}(u,\chi_{D})=\det (I-u\sqrt{q}\Theta_{D})$$
for a unitary $2g\times 2g$ matrix $\Theta_{D}.$
\\By taking a logarithmic derivative of the identity
$$\det (I-u\sqrt{q}\Theta_{D})=(1-u)^{-\lambda}\prod_{P}(1-\chi_{D}(P)u^{\deg P})^{-1},$$
We find
$$-\tr \Theta_{D}^{n}=\frac{\lambda}{q^{n/2}}+\frac{1}{q^{n/2}}\sum_{\deg f=n }\Lambda(f)\chi_{D}(f).\eqno(2.8)$$
\subsection{The Weil bound}
Assume that $B$ is monic of positive degree and not a perfect square. Then we have a bound for the character sum over primes:
$$|\sum_{\substack{\deg P=n\\ P \pr}}(\frac{B}{P})|\ll\frac{\deg B}{n}q^{n/2}.\eqno(2.9)$$
This is deduced from the explicit formula (2.8) when writing $B=DC^{2}$ with $D$ square free of positive degree, and from the unitarity of $\Theta_{D}.$
\section{The hyperelliptic ensemble $\mathcal{H}_{2g+1}$}
\subsection{Averaging over $\mathcal{H}_{2g+1}$}
We denote by $\mathcal{H}_{d}$ the set of square-free monic polynomials of degree $d$ in $\mathbb{F}_{q}[x].$
By using (2.2) and writing
$$\sum_{d\geq 0}\frac{\#\mathcal{H}_{d}}{q^{ds}}=\sum_{f}|f|^{-s}=\frac{\zeta_{q}(s)}{\zeta_{q}(2s)}$$
the sum is over monic and square-free polynomials, we have
$$ \#\mathcal{H}_{d}=\left\{ \begin{array}{ll}
(1-1/q)q^{d} &\mbox{$d\geq2$}\\
q &\mbox{$d=1$}
\end{array}
\right.$$
In particular, for $g\geq 1$,
$$\#\mathcal{H}_{2g+1}=(q-1)q^{2g}.$$

We consider $\mathcal{H}_{2g+1}$ as a probability space with the uniform probability measure, so that the expected value of any function $F$ on $\mathcal{H}_{2g+1}$ is defined as
$$\langle F\rangle:=\frac{1}{\#\mathcal{H}_{2g+1}}\sum_{Q\in \mathcal{H}_{2g+1}}F(Q)\eqno(3.1)$$
We can pick out square-free polynomials by using the M\"{o}bius function $\mu$ of $\mathbb{F}_{q}[x]$
$$\sum_{A^{2}|Q}\mu(A)=
\left\{ \begin{array}{ll}
1 &\mbox{$Q$ is square-free}\\
0 &\mbox{otherwise}
\end{array}
\right.$$
Thus we may write the expected value as
$$\langle F(Q)\rangle=\frac{1}{(q-1)q^{2g}}\sum_{2\alpha+\beta=2g+1}\sum_{\deg B=\beta}\sum_{\deg A=\alpha}\mu(A)F(A^{2}B)\eqno(3.2)$$
the sum is over all monic $A,B.$
\subsection{Averaging quadratic characters}
For a given polynomial $f\in\mathbb{F}_{q}[x]$ apply (3.2) to the quadratic character $\chi_{Q}(f)$. Then
$$\chi_{A^{2}B}(f)=(\frac{B}{f})(\frac{A}{f})^{2}=
\left\{ \begin{array}{ll}
(\frac{B}{f}) &\mbox{$\gcd(A,f)=1$ }\\
0 &\mbox{otherwise}
\end{array}
\right.$$
Hence
$$ \langle \chi_{Q}(f)\rangle=\frac{1}{(q-1)q^{2g}}\sum_{2\alpha+\beta=2g+1}\sum_{\substack{\deg A=\alpha\\ \gcd(A,f)=1}}\mu(A)\sum_{\deg B=\beta}(\frac{B}{f})\eqno(3.3)$$
\subsection{A sum of M$\ddot{\mob}$bius values.}
Define
$$\sigma(f,\alpha):=\sum_{\substack{\deg~A=\alpha \\\gcd(A,f)=1}}\mu(A)\eqno(3.4)$$
Note that $\sigma(f,\alpha)$ depends only on the degrees of the primes dividing $f$, hence we can write for $p_{1},\ldots,p_{n}$ distinct primes of degrees $k_{1},\ldots,k_{n}$ respectively: $\sigma(\prod_{i=1}^{n}p_{i},\alpha)=\sigma(k_{1},\ldots,k_{n};\alpha)$
\begin{lemma}
$$\sigma(k_{1},\ldots,k_{n};\alpha)=\left\{ \begin{array}{ll}
1 &\mbox{$\alpha=0$}\\
-q &\mbox{$\alpha=1$}\\
0 &\mbox{$2\leq \alpha<\min(k_{1},\ldots,k_{n})$}
\end{array}
\right.$$
\end{lemma}
\textbf{Proof:}
By definition
$$\sigma(k_{1},\ldots,k_{n};\alpha)=\sum_{\substack{\deg A=\alpha
\\\gcd(A,p_{1}\cdots p_{n})=1}}\mu(A)$$
Now if $\deg A<\min(k_{1},\ldots,k_{n})$ then $A$ is automatically
coprime to $p_{1},\ldots,p_{n}$ hence in this case the sum is over
all $A$ with degree $\alpha.$ Therefor
$$\sigma(k_{1},\ldots,k_{n};\alpha)=\sum_{\deg A=\alpha}\mu(A)$$
which vanishes if $\alpha \geq 2$ , equals 1 for $\alpha=0$ and $-q$
for $\alpha=1$. $\Box$
\subsection{The probability that $f\nmid Q$}
\begin{lemma}
Let $f=P_{1}P_{2}\cdots P_{k}$ with $P_{1},\ldots,P_{n}$ prime polynomials. Then
$$\langle \chi_{Q}(f^{2})\rangle=1+O(\sum_{P|f}\frac{1}{\|P\|})\eqno(3.5)$$
\end{lemma}
\textbf{Proof:}
\\We may write
$$\chi_{Q}(p_{1}^{2}\cdots p_{k}^{2} )=1-\delta(Q,p_{1}\cdots p_{l} ),~~~~\delta(Q,p_{1}\cdots p_{l} )=\left\{ \begin{array}{ll}
1 &\mbox{$gcd(Q,p_{1}^{2}\cdots p_{l}^{2} )\neq 1$}\\
0 &\mbox{$gcd(Q,p_{1}^{2}\cdots p_{l}^{2} )= 1$}
\end{array}
\right.$$
and hence
$$\langle\chi_{Q}(P_{1}^{2}\cdots P_{n}^{2})\rangle=1-\frac{\#\{Q\in \mathcal{H}_{2g+1}:\exists P_{j}|Q\} }{\#\mathcal{H}_{2g+1}}$$
Replacing the set of square-free $Q$ by arbitrary monic $Q$ of degree $2g + 1$ gives
$$\#\{Q\in \mathcal{H}_{2g+1}:\exists P_{j}|Q\} \leq \#\{\deg Q=2g+1=:\exists P_{j}|Q\}\leq \sum_{j=1}^{k}\frac{q^{2g+1}}{|P_{j}|}$$
so that recalling $\#H_{2g+1}=(q-1)q^{2g}$ we have
$$1-\frac{1}{(1-1/q)}\sum_{j=1}^{k}\frac{1}{|P_{j}|}\geq\langle\chi_{Q}(P_{1}^{2}\cdots P_{n}^{2})\rangle\leq1$$
Thus
$$\langle\chi_{Q}(P_{1}^{2}\cdots P_{n}^{2})\rangle=1+O(\sum_{j=1}^{k}\frac{1}{|P_{j}|})$$
as claimed.
\section{Multiple character sums.}
Define
$$S(\beta;k_{1},\ldots,k_{n}):=
\sum_{\deg p_{1}=k_{1}}\sum_{\deg p_{2}=k_{2}}\ldots\sum_{\deg
p_{n}=k_{n}}\sum_{\deg B=\beta}(\frac{B}{p_{1}p_{2}\cdots p_{n}})\eqno(4.1)$$
the sum over distinct primes
$p_{1},\ldots,p_{n}$ and arbitrary monic $B$.
\\Let $F$ be a squarefree polynomial and
$$\mathcal{L}(u,\chi_{F})=\sum_{\beta=0}^{\infty}A_{F}(\beta)u^{\beta},~~~~A_{F}(\beta):=\sum_{\deg B=\beta}\chi_{F}(B)$$
By quadratic reciprocity (see (2.5))
$$S(\beta;k_{1},\ldots,k_{n})=
(-1)^{\frac{q-1}{2}\beta(\sum_{i=1}^{n}p_{i})}\sum_{\deg
p_{1}=k_{1}}\sum_{\deg p_{2}=k_{2}}\ldots\sum_{\deg
p_{n}=k_{n}}A_{\prod_{i=1}^{n}p_{i}}(\beta)$$ the sum over distinct
primes $p_{1},\ldots,p_{n}.$
\\Since the L-function $\mathcal{L}(u,\chi_{F})$ is a polynomial of degree $\deg F~-~1,$ we have
\begin{lemma}  
If $\beta\geq\sum_{i=1}^{n}k_{i}$ then $S(\beta;k_{1},\ldots,k_{n})=0$
\end{lemma}

\section{Averaging $\prod_{i=1}^{n}(\tr U^{k_{i}})^{a_{i}}$}

\subsection{Reducing to prime powers}

By using the explicit formula we can write
$$\tr U^{k}=\frac{-1}{q^{\frac{k}{2}}}\sum_{\deg f=k}\Lambda(f)\chi_{Q}(f)$$
We separate out the contributions of primes $\mathcal{P}_{k}$,
squares of primes $\Box_{k}$, and higher prime powers
$\mathbb{H}_{k}$ to $\tr U^{k}$:
$$(-trU^{k})^{a}=(\mathcal{P}_{k}+\Box_{k}+\mathbb{H}_{k})^{a}=\sum_{i_{1}+i_{2}+i_{3}=a}{a\choose i_{1},i_{2},i_{3}}(\mathcal{P}_{k})^{i_{1}}(\Box_{k})^{i_{2}}(\mathbb{H}_{k})^{i_{3}}$$
where ${a\choose i_{1},i_{2},i_{3}}=\frac{a!}{i_{1}!i_{2}!i_{3}!}.$
Denote
$$\Box(m,k_{j})=\frac{(\frac{k_{j}}{2})^{m}}{q^{\frac{mk_{j}}{2}}}\sum_{\substack{\deg p_{1},\ldots,p_{m}=\frac{k_{j}}{2}
\\p_{1},\ldots, p_{m}~\dis}}\chi_{Q}((p_{1})^{2}\cdots(p_{m})^{2})\eqno(5.1)$$
$$\Delta_(2m,k_{j})=\frac{(k_{j})^{2m}}{q^{mk_{j}}}\sum_{\substack{\deg p_{1},\ldots,p_{m}=k_{j}
\\p_{1},\ldots, p_{m}~\dis}}\chi_{Q}((p_{1})^{2}\cdots(p_{m})^{2})\eqno(5.2)$$
$$\mathcal{P}(m,k_{j})=\frac{(k_{j})^{m}}{q^{\frac{mk_{j}}{2}}}\sum_{\substack{\deg p_{1},\ldots,p_{m}=k_{j}
\\p_{1},\ldots ,p_{m}~\dis}}\chi_{Q}(p_{1}\cdots p_{m})\eqno(5.3)$$
\\Hence the product  $\prod_{j=1}^{n}(\tr U^{k_{j}})^{a_{j}}$ gives various terms
\\
\\1)Squares
$\prod_{j=1}^{n}\sum_{i_{j}=0}^{\lfloor\frac{a_{j}}{2}\rfloor}{a_{j}\choose
2i_{j}}\frac{(2i_{j})!}{2^{i_{j}}}(a_{j}-2i_{j})!\Delta(2i_{j},k_{j})\Box(a_{j}-2i_{j},k_{j}).$
\\
\\2)Distinct primes $\prod_{j=1}^{n}\mathcal{P}(a_{j},k_{j}).$
\\
\\3)Mixed terms- distinct primes and squares and some high powers: \\$\prod_{j=1}^{n}\sum_{i_{j_{1}}+i_{j_{2}}=a_{j}}(\sum_{m=1}^{i_{j_{1}}}\mathcal{P}(m,k_{j})\Delta(i_{j_{1}}-m,k_{j}))
\Box(i_{j_{2}},k_{j}).$ We can examine a specific term such as
\\$\prod_{j=1}^{n}\mathcal{P}(m_{j},k_{j})\Delta(i_{j_{1}}-m_{j},k_{j})\Box(i_{j_{2}},k_{j})$
since the mixed terms are a finite sum of this kind of terms.
\\
\\4)Higher powers:
\\$\prod_{i=1}^{n}\frac{k_{i}/d_{i_{1}}\cdots
k_{i}/d_{i_{m_{i}}}}{q^{\frac{k_{i}m_{i}}{2}}} \sum_{\deg
p_{i_{j}}=k_{i}}\chi_{Q}(p_{i_{1}}^{d_{i_{1}}}\cdots
p_{i_{m_{i}}}^{d_{i_{m_{i}}}})$ where there is $1\leq i
\leq n$ ; $1\leq j \leq m_{i}$ (at least one index) such that $d_{i_{j}}\geq 3.$
\\
\\Our findings are, assuming $\sum_{j=1}^{n}a_{j}k_{j}\leq 2g-1$ and $\log_{q}g\ll\min(k_{1}a_{1},\ldots,k_{n}a_{n})$:
$$\langle \prod_{j=1}^{n}\sum_{i_{j}=0}^{\lfloor\frac{a_{j}}{2}\rfloor}{a_{j}\choose
2i_{j}}\frac{(2i_{j})!}{2^{i_{j}}}(a_{j}-2i_{j})!\Delta(2i_{j},k_{j})\Box(a_{j}-2i_{j},k_{j})\rangle=$$
$$\prod_{j=1}^{n}\sum_{i_{j}=0}^{\lfloor\frac{a_{j}}{2}\rfloor}(k_{j})^{i_{j}}{a_{j}\choose
2i_{j}}\frac{(2i_{j})!}{2^{i_{j}}(i_{j})!}(-\eta_{k_{j}})^{a_{j}-2i_{j}}
+O(\frac{1}{q^{g}})$$
All the other terms contribute $o(1)$ under these terms.
\subsection{Squares}
Consider the term
$$\prod_{j=1}^{n}\sum_{i_{j}=0}^{\lfloor\frac{a_{j}}{2}\rfloor}{a_{j}\choose
2i_{j}}\frac{(2i_{j})!}{2^{i_{j}}}(a_{j}-2i_{j})!\Delta(2i_{j},k_{j})\Box(a_{j}-2i_{j},k_{j})\eqno(5.4)$$
which equals to
$$\sum_{i_{1}=0}^{\lfloor\frac{a_{1}}{2}\rfloor}\ldots\sum_{i_{n}=0}^{\lfloor\frac{a_{n}}{2}\rfloor}\prod_{j=1}^{n}{a_{j}\choose
2i_{j}}\frac{(2i_{j})!}{2^{i_{j}}}(a_{j}-2i_{j})!\Delta(2i_{j},k_{j})\Box(a_{j}-2i_{j},k_{j}).$$
Hence it is enough to compute the expected value of the term
$\prod_{j=1}^{n}\Delta(2i_{j},k_{j})\Box(a_{j}-2i_{j},k_{j}).$
This contributes to the product
$\prod_{j=1}^{n}(\tr U^{k_{j}})^{a_{j}}$
$$\prod_{j=1}^{n}\frac{(k_{j})^{2i_{j}}(\frac{-k_{j}}{2})^{a_{j}-2i_{j}}}{q^{\frac{k_{j}}{2}a_{j}}}
\sum_{\substack{\deg p_{l,i_{j}}=k_{j}~1\leq l \leq i_{j}\\\deg
p_{l,i_{j}}=\frac{k_{j}}{2}~2i_{j}\leq l \leq a_{j}}}\chi_{Q}(\prod_{j=1}^{n}(p_{1,j}^{2}\cdots p_{i_{j},j}^{2})(p_{2i_{j},j}^{2}\cdots p_{a_{j},j}^{2}))$$
The sum is over different primes.
To average we use Lemma 3.1
$$\langle \chi_{Q}(f^{2})\rangle=1+O(\sum_{P|f}\frac{1}{\|P\|})$$
Hence the average of $\prod_{j=1}^{n}\Delta(2i_{j},k_{j})\Box(a_{j}-2i_{j},k_{j})$ is
$$\prod_{j=1}^{n}\frac{(k_{j})^{2i_{j}}(\frac{-k_{j}}{2})^{a_{j}-2i_{j}}}{q^{\frac{k_{j}}{2}a_{j}}}{\pi(k_{j})\choose i_{j}}{\pi(\frac{k_{j}}{2})\choose a_{j}-2i_{j}}(1+O(\frac{1}{q^{\min{(k_{1},\ldots,k_{n})}}})).$$
The contribution of squares to the product
$\prod_{j=1}^{n}(\tr U^{k_{j}})^{a_{j}}$ is
$$\sum_{i_{1}=0}^{\lfloor\frac{a_{1}}{2}\rfloor}\ldots\sum_{i_{n}=0}^{\lfloor\frac{a_{n}}{2}\rfloor}\prod_{j=1}^{n}\frac{(k_{j})^{2i_{j}}(\frac{-k_{j}}{2})^{a_{j}-2i_{j}}}{q^{\frac{k_{j}}{2}a_{j}}}{a_{j}\choose
2i_{j}}\frac{(2i_{j})!}{2^{i_{j}}}(a_{j}-2i_{j})!{\pi(k_{j})\choose i_{j}}{\pi(\frac{k_{j}}{2})\choose a_{j}-2i_{j}}(1+O(\frac{1}{q^{\min{(k_{1},\ldots,k_{n})}}}))=$$
$$\prod_{j=1}^{n}\sum_{i_{j}=0}^{\lfloor\frac{a_{j}}{2}\rfloor}\frac{(k_{j})^{2i_{j}}(\frac{-k_{j}}{2})^{a_{j}-2i_{j}}}{q^{\frac{k_{j}}{2}a_{j}}}{a_{j}\choose
2i_{j}}\frac{(2i_{j})!}{2^{i_{j}}}(a_{j}-2i_{j})!{\pi(k_{j})\choose i_{j}}{\pi(\frac{k_{j}}{2})\choose a_{j}-2i_{j}}(1+O(\frac{1}{q^{\min{(k_{1},\ldots,k_{n})}}}))=$$
$$\prod_{j=1}^{n}\sum_{i_{j}=0}^{\lfloor\frac{a_{j}}{2}\rfloor}\frac{(k_{j})^{2i_{j}}(\frac{-k_{j}}{2})^{a_{j}-2i_{j}}}{q^{\frac{k_{j}}{2}a_{j}}}{a_{j}\choose
2i_{j}}\frac{(2i_{j})!}{2^{i_{j}}}\frac{\pi(k_{j})!\pi(\frac{k_{j}}{2})!}{(i_{j})!(\pi(k_{j})-i_{j})!(\pi(\frac{k_{j}}{2})-a_{j}+2i_{j})!}(1+
O(\frac{1}{q^{\min{(k_{1},\ldots,k_{n})}}}))=$$
$$\prod_{j=1}^{n}\sum_{i_{j}=0}^{\lfloor\frac{a_{j}}{2}\rfloor}(k_{j})^{i_{j}}{a_{j}\choose
2i_{j}}\frac{(2i_{j})!}{2^{i_{j}}(i_{j})!}(-\eta_{k_{j}})^{a_{j}-2i_{j}}
+o(1)\eqno(5.5)$$
\subsection{Primes}
In this section we focus on the contribution of different primes:
Notice that the case of different primes is equivalent to the case of $a_{j}=1,1\leq j\leq n$ and all the $k_{j}'s$ are different. Hence we consider the case of $\prod_{i=1}^{n}\mathcal{P}_{k_{i}}.$
Assume (for the convenience of writing)  $k_{1}=\min (k_{1},\ldots,k_{n})$ .
We use (3.2) and the explicit formula of (2.8) for the
mean value of $\prod_{i=1}^{n}\mathcal{P}_{k_{i}}$:
$$\langle\prod_{i=1}^{n}\mathcal{P}_{k_{i}}\rangle=\frac{(-1)^{n}(\prod_{i=1}^{n}k_{i})}{q^{\frac{\sum_{i=1}^{n}k_{1}}{2}+2g}(q-1)}
\sum_{\deg P_{i}=k_{i}}\sum_{2\alpha+\beta=2g+1}\sum_{\substack{\deg
A=\alpha\\ \gcd(A,P_{i})=1}}\mu(A)\sum_{\deg
B=\beta}(\frac{B}{\prod_{i=1}^{n}P_{i}})=$$
$$\frac{(-1)^{n}(\prod_{i=1}^{n}k_{i})}{q^{\frac{\sum_{i=1}^{n}k_{1}}{2}+2g}(q-1)}
\sum_{0\leq\alpha\leq
g}\sigma(k_{1},\ldots,k_{n};\alpha)S(2g+1-2\alpha;k_{1},\ldots,k_{n})\eqno(5.6)$$
Provided $\sum_{i=1}^{n}k_{i}<2g$, the case $\alpha=0$ and the case
$\alpha=1$ give by Lemma 4.1
$S(2g+1-2\alpha;k_{1},\ldots,k_{n})=0$. Now
$\sigma(k_{1},\ldots,k_{n};\alpha)=0$ for $2\leq\alpha\leq k_{1}$ by Lemma 3.1. Thus it suffices to take $k_{1}\leq\alpha$ and
$2g+1-2\alpha<\sum_{i=1}^{n}k_{i}$.
\\For this we use the Weil bound
$$S(2g+1-2\alpha;k_{1},\ldots,k_{n})\ll\frac{(2g+1-2\alpha)^{n}}{k_{1}k_{2}\cdots k_{n}}q^{2g+1-2\alpha+\frac{\sum_{i=1}^{n}k_{i}}{2}}\eqno(5.7)$$
To get
$$\langle\prod_{i=1}^{n}\mathcal{P}_{k_{i}}\rangle=\frac{(\prod_{i=1}^{n}k_{i})}{q^{\frac{\sum_{i=1}^{n}k_{1}}{2}+2g}(q-1)}
\sum_{k_{1}\leq\alpha\leq
g}\sigma(k_{1},\ldots,k_{n};\alpha)S(2g+1-2\alpha;k_{1},\ldots,k_{n})$$
$$\ll\frac{q(-1)^{n}}{(q-1)}
\sum_{max(k_{1},g-\frac{\sum_{i=1}^{n}k_{i}-1}{2})\leq\alpha\leq
g}\sigma(k_{1},\ldots,k_{n};\alpha)q^{-2\alpha}(2g+1-2\alpha)^{n}\eqno(5.8)$$
Notice that $\sigma(k_{1},\ldots,k_{n};\alpha)\ll q^{\alpha}$,hence the  above term is bounded by
$$\frac{g^{n+1}}{q^{max(k_{1},g-\frac{\sum_{i=1}^{n}k_{i}-1}{2})-1}}$$
Provided
$(n+2)\log_{q}g < k_{1}$ this is $o(1).$
\subsection{Mixed terms: Primes and Squares}
Define $\prod_{j=1}^{n}\mathcal{P}(i_{j}-2m_{j},k_{j})\Delta(2m_{j},k_{j})\Box(a_{j}-i_{j},k_{j})$
to be the contribution of primes, squares and some higher powers, to $\prod_{j=1}^{n}(\tr U^{k_{j}})^{a_{j}}.$ In this case $i_{j}-2m_{j}\neq 0$ for at least one of the $j's$ , and for $j$ such that $k_{j}$ is odd we have $i_{j}=a_{j}.$ For the convenience of writing we will bound the expected value of the term $\mathcal{P}(i-2m,k)\Delta(2m,k)\Box(a-i,k).$
The expected value of the product will be bounded exactly in the same methods.
We start by writing
$$\chi_{Q}(p_{1}^{2}\cdots p_{l}^{2} )=1-\delta(Q,p_{1}\cdots p_{l} ),~~~~\delta(Q,p_{1}\cdots p_{l} )=\left\{ \begin{array}{ll}
1 &\mbox{$gcd(Q,p_{1}^{2}\cdots p_{l}^{2} )\neq 1$}\\
0 &\mbox{$gcd(Q,p_{1}^{2}\cdots p_{l}^{2} )= 1$}
\end{array}
\right.$$ Define
$$\omega_{Q}^{l}(k_{j}):=\sum_{\substack{\deg p_{1}\cdots p_{l}=k \\p_{1}\ldots p_{l}~\dis}}\delta(Q,p_{1}\cdots p_{l} )\eqno(5.9)$$
The sum is over prime factors. This satisfies
$$\omega_{Q}^{l}(k_{j})\leq \frac{2g+1}{k_{j}}\pi(k_{j})^{l-1}$$

$$\Box(a-i,k)=\frac{(\frac{k}{2})^{a-i}}{q^{(a-i)\frac{k}{2}}}\sum_{\substack{\deg p_{1},\ldots,p_{a-i}=\frac{k}{2}
\\p_{1},\ldots,p_{a-i}~\dis}}\chi_{Q}((p_{1})^{2}\cdots(p_{a-i})^{2})=
\frac{(\frac{k}{2})^{a-i}}{q^{(a-i)\frac{k}{2}}}({\pi(\frac{k}{2})\choose a-i}-\omega_{Q}^{a-i}(\frac{k}{2}))$$

$$\Delta(2m,k)=\frac{k^{2m}}{q^{mk}}\sum_{\substack{\deg p_{1},\ldots,p_{m}=k
\\p_{1},\ldots,p_{m}~\dis}}\chi_{Q}((p_{1})^{2}\cdots(p_{m})^{2})=
\frac{k^{2m}}{q^{mk}}({\pi(k)\choose m}-\omega_{Q}^{m}(k))$$
Hence
$$\mathcal{P}(i-2m,k)\Delta(2m,k)\Box(a-i,k)=
\mathcal{P}(i-2m,k)
\frac{k^{2m}}{q^{mk}}({\pi(k)\choose m}-\omega_{Q}^{m}(k))
\frac{(\frac{k}{2})^{a-i}}{q^{(a-i)\frac{k}{2}}}({\pi(\frac{k}{2})\choose a-i}-\omega_{Q}^{a-i}(\frac{k}{2}))$$
Notice that
$$\frac{k^{2m}}{q^{mk}}{\pi(k)\choose m}=\frac{k^{m}}{m!}+O(\sum_{l=1}^{m-1}\frac{k^{m+l}}{q^{lk}})\eqno(5.10)$$
$$\frac{(\frac{k}{2})^{a-i}}{q^{(a-i)\frac{k}{2}}}{\pi(\frac{k}{2})\choose a-i}=
\frac{1}{(a-i)!}+O(\sum_{l=1}^{a-i-1}\frac{\frac{k}{2}^{l}}{q^{l\frac{k}{2}}})\eqno(5.11)$$
By (5.10),(5.11) and the prime section we have
$$\langle\mathcal{P}(i-2m,k)
\frac{k^{2m}}{q^{mk}}{\pi(k)\choose m}
\frac{(\frac{k}{2})^{a-i}}{q^{(a-i)\frac{k}{2}}}{\pi(\frac{k}{2})\choose a-i}\rangle\sim
\frac{k^{m}}{m!(a-i)!}\langle\mathcal{P}(i-2m,k)\rangle\ll \frac{k^{m}}{m!(a-i)!}\frac{g^{i-2m+1}}{q^{\max(k,g-(i-2m)\frac{k}{2})}}$$
which is $o(1)$ provided $k>a\log_{q}g.$
\\It is enough to compute the expected value of
$$\langle\mathcal{P}(i-2m,k)
\frac{k^{2m}}{q^{mk}}\omega_{Q}^{m}(k)
\frac{(\frac{k}{2})^{a-i}}{q^{(a-i)\frac{k}{2}}}\omega_{Q}^{a-i}(\frac{k}{2})\rangle$$
By the Cauchy-Schwartz inequality
$$\langle\mathcal{P}(i-2m,k)
\frac{k^{2m}}{q^{mk}}\omega_{Q}^{m}(k)
\frac{(\frac{k}{2})^{a-i}}{q^{(a-i)\frac{k}{2}}}\omega_{Q}^{a-i}(\frac{k}{2})\rangle\leq
\langle(\mathcal{P}(i-2m,k))^{2}\rangle\frac{(2g+1)^{2}k^{m}}{q^{\frac{3k}{2}}}\eqno(5.12)$$
next we show that $\langle(\mathcal{P}(i-2m,k))^{2}\rangle$ is polynomial in $g.$ It follow that for $ \log_{q}g\ll k$ the above is $o(1)$.
$$\langle(\mathcal{P}(i-2m,k))^{2}\rangle=
\langle((\mathcal{P}_{k})^{(i-2m)}-\Delta(2,k)(\mathcal{P}_{k})^{(i-2m-2)})^{2}\rangle=$$
$$\langle((\mathcal{P}_{k})^{(2i-4m)}-2\Delta(2,k)(\mathcal{P}_{k})^{2i-4m-2}+(\Delta(2,k))^{2}(\mathcal{P}_{k})^{2(i-2m-2)})\rangle$$
We use on this term the same methods as before to have
$$\langle((\mathcal{P}_{k})^{(2i-4m)}\rangle-
2\langle(\frac{k^{2}}{q^{k}}(\pi(k)-\omega_{Q}^{1}(k)))(\mathcal{P}_{k})^{2i-4m-2}\rangle+
\langle((\frac{k^{2}}{q^{k}}(\pi(k)-\omega_{Q}^{1}(k)))^{2}(\mathcal{P}_{k})^{2(i-2m-2)})\rangle$$
This comes down to bounding the general term $\langle(\mathcal{P}_{k})^{2l}\rangle,$ since
$$\langle\frac{k^{2}}{q^{k}}\pi(k)(\mathcal{P}_{k})^{2i-4m-2}\rangle\sim k\langle(\mathcal{P}_{k})^{2i-4m-2}\rangle~~;~~
\langle(\frac{k^{2}}{q^{k}}\pi(k))^{2}(\mathcal{P}_{k})^{2(i-2m-2)})\rangle\sim k^{2}\langle(\mathcal{P}_{k})^{2(i-2m-2)})\rangle$$
and for terms with $\omega_{Q}^{1}(k)$ we use the Cauchy-Schwartz inequality. For example
$$\langle\frac{k^{2}}{q^{k}}\omega_{Q}^{1}(k)(\mathcal{P}_{k})^{2i-4m-2}\rangle\leq \frac{(2g+1)k}{q^{k}}\langle(\mathcal{P}_{k})^{4i-8m-4}\rangle^{\frac{1}{2}}$$
\begin{lemma} For $2g>k>2l\log_{q}g$
$$\langle(\mathcal{P}_{k})^{2l}\rangle=O(g^{l})\eqno(5.13)$$
\end{lemma}
\textbf{Proof:} We will prove the lemma by induction .
\\For $\textbf{\underline{l=1}}$ we have
$$\langle(\mathcal{P}_{k})^{2}\rangle=\langle\mathcal{P}(2,k)\rangle+
\langle\Delta(2,k)\rangle$$
By section $5.3$ (the prime section) we have
$\langle\mathcal{P}(2,k)\rangle= o(1)$
For the second term we use $(3.5)$
$$\langle\Delta(2,k)\rangle=\frac{k^{2}}{q^{k}}\pi(k)(1+O(\frac{1}{q^{k}}))=k+O(\frac{k}{q^{k}}).$$
In conclusion, for $2g>k>2\log_{q}g$  we have
$$\langle(\mathcal{P}_{k})^{2}\rangle=O(g)$$
For $\textbf{\underline{l=2}}$
$$\langle(\mathcal{P}_{k})^{4}\rangle=\langle\mathcal{P}(4,k)\rangle+
\langle\Delta(2,k)(\mathcal{P}_{k})^{2}\rangle$$
By section $5.3$ (the prime section) we have
$\langle\mathcal{P}(4,k)\rangle=o(1)$
For the second term we use Cauchy-Schwartz inequality
$$\langle\Delta(2,k)(\mathcal{P}_{k})^{2}\rangle=
\langle\frac{k^{2}}{q^{k}}(\pi(k)-\omega_{Q}^{1}(k)))(\mathcal{P}_{k})^{2}\rangle\sim
k\langle(\mathcal{P}_{k})^{2}\rangle-\langle\frac{k^{2}}{q^{k}}\omega_{Q}^{1}(k)(\mathcal{P}_{k})^{2}\rangle\leq$$
$$k\langle(\mathcal{P}_{k})^{2}\rangle+\frac{k(2g+1)}{q^{k}}\langle(\mathcal{P}_{k})^{2}\rangle^{\frac{1}{2}}$$
By the case of $\textbf{\underline{l=1}}$ we have
$$\langle(\mathcal{P}(4,k)\rangle =O(g^{2})$$
For the case of \underline{general $l$}:
$$\langle(\mathcal{P}_{k})^{2l}\rangle=\langle\mathcal{P}(2l,k)\rangle+
\langle\Delta(2,k)(\mathcal{P}_{k})^{2l-2}\rangle$$
By section $5.3$ (the prime section) we have
$\langle\mathcal{P}(2l,k)\rangle=o(1).$
For the second term we use the same method as before to get
$$\langle\Delta(2,k)(\mathcal{P}_{k})^{2l-2}\rangle\leq k\langle(\mathcal{P}_{k})^{2l-2}\rangle+\frac{k(2g+1)}{q^{k}}\langle(\mathcal{P}_{k})^{2l-2}\rangle$$
By the induction
$$\langle(\mathcal{P}_{k})^{2l-2}\rangle=O(g^{l-1})$$
Hence, provided $2l\log_{q}g<k<2g$ we have
$$\langle(\mathcal{P}_{k})^{2l}\rangle=O(g^{l})$$
this concludes the lemma.
\\
\\
\\Going back we get $\langle(\mathcal{P}(i-2m,k))^{2}\rangle$ is polynomial in $g.$ It follows that the contribution of mixed terms of primes and squares to the expected value of the product of traces is $o(1).$
\subsection{Higher powers}
We now consider the contribution of higher powers to the product of traces. These arise from terms with $\mathbb{H}_{k_{j}}$ or from the remaining terms which were not considered before (notice that some of the cases involving higher powers were considered in the previous section).
These terms coincides up to division by factors such as $q^{\frac{k_{j}}{2}(1-\frac{1}{d})}$ ($d$ is some finite integer), with the contribution of one of the previous types. Since the mean value of these previous terms is in all cases bounded by $g^{\sum_{i=1}^{n}a_{i}}$, after the division we get a negligible term (provided $\log_{q}g\ll\min(k_{1},\ldots,k_{n})$).
\subsection{Conclusion of the proof}
We saw that $\langle\prod_{j=1}^{n}(\tr U^{k_{j}})^{a_{j}}\rangle$ is the sum of the expected values of the Frobenius class of hyperelliptic curves of genus $g$ over the field $\mathbb{F}_{q}$ of various terms:
\\
\\1)squares
$\prod_{j=1}^{n}\sum_{i_{j}=0}^{\lfloor\frac{a_{j}}{2}\rfloor}{a_{j}\choose
2i_{j}}\frac{(2i_{j})!}{2^{i_{j}}}(a_{j}-2i_{j})!\Delta(2i_{j},k_{j})\Box(a_{j}-2i_{j},k_{j}).$
\\
\\2)distinct primes $\prod_{j=1}^{n}\mathcal{P}(a_{j},k_{j}).$
\\
\\3)mixed terms- distinct primes and squares and some high powers. \\$\prod_{j=1}^{n}\sum_{i_{j_{1}}+i_{j_{2}}=a_{j}}(\sum_{m=1}^{i_{j_{1}}}\mathcal{P}(m,k_{j})\Delta(i_{j_{1}}-m,k_{j}))
\Box(i_{j_{2}},k_{j}).$
\\
\\4)$\mathbb{H}_{k_{1},\ldots,k_{n}}$ Higher powers.
\\
\\We saw that under the following condition: $k_{j}\in \{1,2,\ldots\}$ for $1\leq j \leq n$ are such that $\sum_{j=1}^{n}a_{j}k_{j}\leq 2g-1$ for fixed integers $a_{j}$, and $\log_{q}g\ll\min(k_{1},\ldots,k_{n})$,  the expected value of all the above terms is negligible - $o(1)$, except from the first term, the squares. This gives
$$\langle \prod_{j=1}^{n}\sum_{i_{j}=0}^{\lfloor\frac{a_{j}}{2}\rfloor}{a_{j}\choose
2i_{j}}\frac{(2i_{j})!}{2^{i_{j}}}(a_{j}-2i_{j})!\Delta(2i_{j},k_{j})\Box(a_{j}-2i_{j},k_{j})\rangle=$$
$$\prod_{j=1}^{n}\sum_{i_{j}=0}^{\lfloor\frac{a_{j}}{2}\rfloor}(k_{j})^{i_{j}}{a_{j}\choose
2i_{j}}\frac{(2i_{j})!}{2^{i_{j}}(i_{j})!}(-\eta_{k_{j}})^{a_{j}-2i_{j}}
+O(\frac{1}{q^{g}})$$
Putting this together gives Theorem 1.1.

\end{document}